%% file: main.tex
\documentclass[preprint,12pt]{elsarticle}
\usepackage{amssymb}
\usepackage{amsmath}

\usepackage{amsfonts}  
\usepackage{physics}   
\usepackage{parskip}
\usepackage{hyperref}

\hypersetup{
    colorlinks=true,
    linkcolor=blue,
    filecolor=magenta,
    urlcolor=cyan,
    pdftitle={Time Marching for Couple Stress Elastodynamics},
    pdfauthor={José Ortiz-Ocampo, Nicolás Guarín-Zapata},
    pdfsubject={Consistent couple-stress theory, Dynamic problems, Micromechanics, Mixed Finite Element Methods}
}

\journal{Computer Methods in Applied Mechanics and Engineering}

\begin{document}

\begin{frontmatter}



\title{\textbf{Implicit Time-Marching for Lagrange Multiplier Formulation for Couple Stress Elastodynamics}}


\author[EAFIT]{José Ortiz-Ocampo}
\ead{jhortizo@eafit.edu.co}
\author[EAFIT]{Nicolás Guarín-Zapata}
\ead{nguarinz@eafit.edu.co}

\affiliation[EAFIT]{organization={School of Applied Sciences and Engineering, Universidad EAFIT},
            addressline={Cra. 48 \#4 Sur-414}, 
            city={Medellín},
            postcode={050022}, 
            state={Antioquia},
            country={Colombia}}

\begin{abstract}
The study of metamaterials and architected materials has intensified interest in continuum mechanics models that capture size-dependent microstructure interactions. Among these, Consistent Couple-Stress Theory (C-CST) incorporates microscale mechanical interactions by introducing higher-order derivatives in the strain energy. While previous studies have relied on convolutional principles or inverse Laplace transforms to obtain time-dependent solutions, this work demonstrates that implicit time integration applied to a mixed finite element method with a Lagrange multiplier provides stable, direct time-domain solutions for dynamic C-CST modeling. The proposed finite element scheme is tested through the Method of Manufactured Solutions (MMS) for static cases and dynamic simulations of simple mechanical scenarios.

Our computational experiments revealed energy dissipation, emphasizing the importance of exploring symplectic integrators in future work to impose energy conservation. Additionally, further research is required to verify temporal stability through time-domain MMS and to investigate complex mechanical scenarios, including those previously restrictive, challenging to simulate, or unfeasible with existing dynamic methods. This work lays the groundwork for studying size-dependent material behavior and provides the foundation for advanced applications in material design and wave propagation.

\end{abstract}




\end{frontmatter}

\input{sections/01_introduction}
\input{sections/02_model}

\input{sections/03_results}
\input{sections/04_conclusions}
\input{sections/appendix}

\section*{Aknowledgement}

José Ortiz-Ocampo thanks Universidad EAFIT for supporting this work through the Master's Scholarship Program \textbf{Honor Posgrado - Maestría}.

\bibliographystyle{plainnat}
\bibliography{references.bib}

\end{document}

%% file: sections/01_introduction.tex
\section{Introduction}
\label{sec:introduction}

The expanding research in metamaterials and architected materials has po\-pularized the study of continuum mechanics models that include local microstructure interactions \cite{nowacki_theory_1971, do_rosario_combined_2020, ulloa_fracture_2024, ariza_homogenization_2024}. Several modern continuum theories address these size-dependent, non-classical phenomena, including couple-stress theories (CSTs) \cite{kumar_size-dependent_2023, mohammadimehr_magneto-mechanical_2021, abouelregal_modified_2024, rahmani_eringens_2022}.

The family of couple-stress models was first proposed in the early 20th century by the Cosserat brothers, incorporating local microrotations as additional degrees of freedom in material points~\cite{cosserat_theorie_1909}. This generalizes Cauchy's postulate by including couple stresses, representing mechanical interactions related to couples per unit surface. These effects become significant when system dimensions are comparable to the material’s microstructural features \cite{dargush_convolved_2023, guarin-zapata_variational_2021, mase_continuum_2011}.

Several extensions and formalizations of the theory emerged in later works, such as the micropolar, microstretch, and micromorphic theories \cite{eringen_linear_1966, nowacki_theory_1971, mindlin_micro-structure_1964, eringen_nonlinear_1964}. At the same time, works from Toupin, Mindlin and others introduced kinematic enrichment by equating local microrotations to macro-rotations, leading to what is known as indeterminate couple-stress theory due to the inability to determine the spherical part of the couple-stress tensor and inconsistencies in the boundary conditions~ like \cite{toupin_elastic_1962, mindlin_effects_1962, koiter_couple-stresses_1969, eringen_theory_1999}.

Hadjesfandiari and Dargush introduced a novel theory that addressed those inconsistencies in the boundary conditions and indetermination of the spherical component~\cite{hadjesfandiari_couple_2011}. This model, known as Consistent Couple-Stress Theory (C-CST), determines all components of the force and couple stress tensors, identifies a necessary and sufficient set of boundary conditions, and eliminates redundant force components.

In recent years, there has been a growing interest in the dynamic behavior predicted by C-CST. \citet{deng_mixed_2016, deng_mixed_2017} investigated natural frequency responses for size-dependent couple-stress using a mixed Lagrange multiplier formulation. Later, Guarín-Zapata et al developed a frequency action principle and implemented a frequency-domain Finite Element Method (FEM) to study the dispersion of waves in phononic crystals~\cite{guarin-zapata_variational_2021}. Dargush introduced a convolved action principle for C-CST elastodynamics and verified it numerically using a Laplace-domain Boundary Element Method (BEM) analysis~\cite{dargush_convolved_2023}. Lei et al provided Laplace-domain BEM numerical solutions and later derived time-dependent solutions using Durbin's method for the inverse Laplace transform~\cite{lei_laplace-domain_2023, durbin_numerical_1974}.

Time integration schemes have been, and remain, widely studied and implemented in elastodynamics \cite{kadapa_novel_2021, rossi_implicit_2016, soares_novel_2022, soares_truly-explicit_2023, song_high-order_2024}. However, to the best of the authors' knowledge, they have not yet been used for C-CST formulation. All the  previously mentioned approaches require workarounds to obtain the time-depen\-dent solution, either through convolutional principles or inverse transformations. Here, we compute the time-dependent solution directly using implicit time integration, which enables iterative and controlled calculation of the dynamic solution. Their implementation on top of the current mixed FEM with Lagrange multiplier method presents additional challenges compared to other elastodynamic models due to the nature and structure of the mass and stiffness matrices. One of the main contributions of this paper is the derivation of equations suitable for implementing time integration, as well as its verification through the analysis of the time evolution of several case studies.

This paper is organized as follows: Section 2 provides an overview of the C-CST model and derives a weak formulation of the problem based on D'Alembert's principle. We present the Finite Element formulation and the derivation of the implicit time integration scheme. Section 3 presents computational tests, including static analysis to check the method and matrices, as well as dynamic cases to assess the scheme. Section 4 concludes the paper.

%% file: sections/02_model.tex

\section{Mathematical modeling}
\label{sec:mathematical_modeling}

In this section, we present the mathematical procedures to derive the C-CST time-domain equations based on the D'Alembert principle, as well as their space and time discretizations via FEM and implicit time integration \cite{guarin-zapata_variational_2021, darrall_finite_2014}.

\subsection{C-CST time-domain equations}
\label{subsec:ccst_equations}

We present a derivation of the time-domain C-CST equations, for a complete description of the C-CST model please refer to
\cite{hadjesfandiari_couple_2011, guarin-zapata_variational_2021}  and the references therein.
Here, we use index notation in cartesian coordinates and assume the summation convention.

We start by considering a continuum body with volume \(V\) and boundary \(S\) with \(n_i\) the vector normal to the surface, in which known body forces \(f_i\), force tractions \(t_i\) and couple-tractions \(m_i\) are applied, where the subindex indicates cartesian components of the vectors. We are not considering body couples since \citet{hadjesfandiari_couple_2011} showed that
they can be written in terms of couple-tractions.
\begin{figure}[ht]
\centering
\includegraphics[width=3 in]{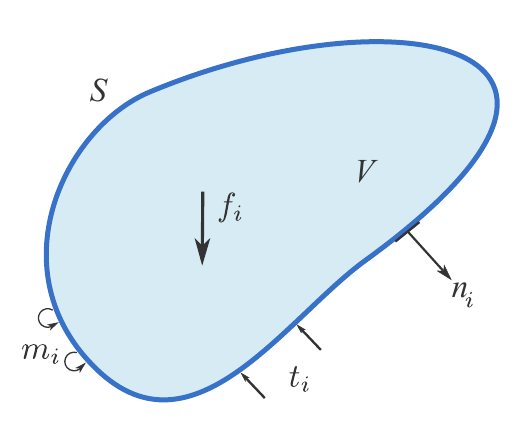}
\caption{Schematic of a CCST body, subject to body forces \(f_i\) applied on the volume, and force and couple tractions \(t_i\) and \(m_i\) applied on the boundary. \(n_i\) indicates the normal vector to the surface.}

\end{figure}

We consider that the continuum has translational and rotational mechanical 
interactions between material points in terms of force-traction 
\(t_i^{(\hat{n})}\) and couple-traction \(m_i^{(\hat{n})}\) vectors defined 
through a generalized Cauchy's postulate \cite{mase_continuum_2011}. Force-tractions and couple-tractions are described by the non-symmetric force-stress \(\sigma_{ij}\) and couple-stress \(\mu_{ij}\) tensors, written as
\begin{align*}
    &t_i^{(\hat{n})} = \sigma_{ji} n_j\, ,\\
    &m_i^{(\hat{n})} = \mu_{ji} n_j\, .
\end{align*}

The conservation of linear and angular momentum is given by
\begin{align*}
    \int\limits_S t_i \dd{S} + \int\limits_V f_i \dd{V} = \int\limits_V \rho \ddot{u}_i \dd{V}\, ,\\
    \int\limits_S (\epsilon_{ijk} x_j t_k + m_i) \dd{S} +
    \int\limits_V \epsilon_{ijk} x_j f_k \dd{V} = \int\limits_V \rho \epsilon_{ijk} x_j \ddot{u}_k \dd{V} \, ,
\end{align*}
with \(x_i\) the position vector, \(u_i\) the displacement field, \(\epsilon_{ijk}\) the Levi-Civita permutation symbol, $\rho$ the continuum density, and using dot notation for the time derivative. Using the definitions for the stresses and applying the divergence theorem, we get
\begin{align*}
    \int\limits_V(\sigma_{ji,j} + f_i - \rho \ddot{u}_i) \dd{V} = 0\, ,\\
    \int\limits_V(\epsilon_{ijk} \sigma_{jk} + \mu_{ji,j}) \dd{V} = 0\, .
\end{align*}

Since the volume \(V\) was arbitrary \cite{nowacki_theory_1971}, we ended up with the following equations.
\begin{align}
    \sigma_{ji,j} + f_i = \rho \ddot{u}_i\, ,\\
    \epsilon_{ijk} \sigma_{jk} + \mu_{ji,j}  = 0 \ ,
    \label{eq:conservation_eqs}
\end{align}
where the first equation presents the conservation of momentum in differential form. The second relates the skew-symmetric part of the stress tensor with the divergence of the couple-stress tensor.

Let us introduce \(e_{ij}\) as the classical infinitesimal strain tensor, \(\theta_{i}\) the rotation field and \(\kappa_{i}\) the mean curvature field, which are given by:
\begin{subequations}
\label{eq:kinematics}
\begin{align}
    e_{ij} &= \frac{1}{2}(u_{i,j} + u_{j,i})\, , \label{eq:strain}\\
    \theta_{i} &= \frac{1}{2}\epsilon_{ijk}u_{k ,j}\, , \label{eq:rotation} \\
    \kappa_{i} &= \frac{1}{2}\epsilon_{ijk}\theta_{k ,j}\, . \label{eq:curvature}
\end{align}
\end{subequations}

In the C-CST model, \(\mu_{ij}\) is skew-symmetric, so it can be also written as
\[\mu_i = \frac{1}{2}\epsilon_{ijk}\mu_{kj}\, .\] Moreover, the constitutive equations for linear, elastic centrosymmetric media are given by
\begin{equation}
    \begin{split}
    &\sigma_{ij} = C_{ijkl} e_{kl}\, ,\\
    &\mu_i = D_{ij} \kappa_j\, ,
    \end{split}
    \label{eq:constitutive_aniso}
\end{equation}
Here \(C_{ijkl}\) is the classical stiffness tensor and \(D_{ij}\) is an additional material tensor that appears in C-CST to account for couple-stress effects. In the isotropic case, they are given by
\begin{equation}
    \begin{split}
    &C_{ijkl} = \lambda \delta_{ij} \delta_{kl} + \mu ( \delta_{ik} \delta_{jl}
       + \delta_{il} \delta_{jk} ), \\
    &D_{ij} = 4\eta \delta_{ij},
    \end{split}
    \label{eq:isotropic_tensors}
\end{equation}
where \(\lambda\), \(\mu\) are the  classical Lamé parameters, and \(\eta\) is a non-classical material parameter related to couple-stress effects. In this model, another parameter of interest appears, given by 
\begin{equation}
    l^2 = \frac{\eta}{\mu}\, ,
    \label{eq:lengthscale}
\end{equation}
which is the intrinsic material length scale. When this parameter is comparable
with the length scale of the geometry, one has couple-stress effects in the response of the solid.

Additionally, the coupling of the conservation equations~(\ref{eq:conservation_eqs}) with the constitutive relations in the case of isotropic materials~(\ref{eq:isotropic_tensors}) allows us to obtain the following differential equation in terms of displacements
\begin{equation}
    \left(\lambda + \mu + \eta\nabla^2\right)u_{k, ki} + \left(\mu - \eta\nabla^2\right)\nabla^2 u_i + f_i = \rho \pdv[2]{u_i}{t} ,
    \label{eq:ccst-pde_index}
\end{equation}

where \(t\) denotes times, in agreement with standard notation. However, it must not be confused with the force-traction \(t_i\) previously defined and distinguished by its subscript and vector nature.

Finally, we present the differential equation in vector notation and their corresponding boundary conditions for the problem to be well-posed~\cite{guarin-zapata_variational_2021, darrall_finite_2014}. The boundary \(S\) is split into different boundary conditions: \(S_u\) represents the region of \(S\) with prescribed displacements, \(S_t\) represents the region with prescribed tractions, \(S_\theta\) represents the region with prescribed rotations, and \(S_m\) the region with prescribed couple-tractions. Additionally, \(S = S_u \cup S_t = S_\theta \cup S_m\) and \(S_u \cap S_t = S_\theta \cap S_m = \emptyset\). In general, \(S_u\) and \(S_t\) might overlap with \(S_\theta\) and \(S_m\).
\begin{figure}[ht]
\centering
\includegraphics[width=5 in]{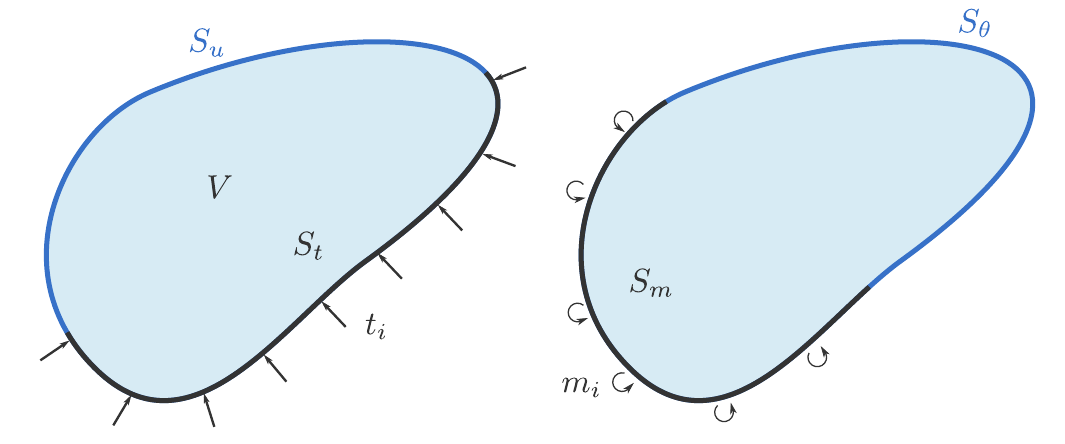}
\caption{Schematic representation of the domain and boundary conditions for the C-CST model. Known body forces \(f_i\) are applied on the volume, and the boundary \(S\) is split into boundaries where different conditions are prescribed.}
\end{figure}

Let us find \(\vb{u}: V \times (0, T] \rightarrow \mathbb{R}\) satisfying the following differential equation
\begin{equation}\label{eq:ccst-pde}
(\lambda + 2\mu)\nabla (\nabla \cdot \vb{u}) - (\mu - \eta \nabla^2)\nabla\times \nabla 
\times \vb{u} + \vb{f} = \rho \ddot{\vb{u}} \, \quad\forall \vb{x} \in V,\, \forall t \in (0, T],
\end{equation}
subject to the following conditions
\begin{align*}
& \vb{u} = \bar{\vb{u}}, &\forall \vb{x} \in S_u\, ,\forall t \in (0, T]\\
&\hat{\vb{n}} \cdot \boldsymbol{\sigma} = \bar{\vb{t}}, &\forall \vb{x} \in S_t\, ,\forall t \in (0, T] \\
&\frac{1}{2}\nabla\times\vb{u} = \bar{\boldsymbol{\theta}}, &\forall \vb{x} \in S_\theta\, ,\forall t \in (0, T]\\
& \boldsymbol{\mu} \cross \hat{\vb{n}} = \bar{\vb{m}}, &\forall \vb{x} \in S_m\, ,\forall t \in (0, T]\\
& \vb{u}(\vb{x}, 0) = \vb{u}_0, &\forall \vb{x} \in V\\
& \dot{\vb{u}}(\vb{x}, 0) = \dot{\vb{u}}_0, &\forall \vb{x} \in V
\end{align*}
where \(\bar{\vb{u}}, \bar{\vb{t}}, \bar{\boldsymbol{\theta}}, \bar{\vb{m}}, \vb{u}_0, \dot{\vb{u}}_0\) are known functions. Notice that, according to equation~(\ref{eq:ccst-pde}), the functions should live in a suitable space of \(u_i \in C^4(V) \times C^2((0, T])\), due to the term \(\nabla^2\, \nabla \times \nabla \times \mathbf{u}\), although the requirement may be less restrictive when analyzed component-wise.

\subsection{Finite element method formulation}
\label{subsec:fem_formulation}

To derive the finite element formulation for C-CST, we begin with the application of the D'Alembert principle over the continuum body \cite{lanczos_variational_1986}. This principle, which generalizes the principle of virtual work to dynamic systems \cite{malvern_introduction_1969, reddy_energy_2017}, states that the total virtual work of all forces in the system, including inertial forces, must vanish for any admissible virtual displacement \(\delta \mathbf{u}\). Specifically, the principle introduces the inertial forces into the virtual work expression, ensuring dynamic equilibrium by requiring that their contributions balance with those of the internal and external forces.

With the definitions presented in Section~\ref{subsec:ccst_equations}, we express D'Alembert principle for the C-CST continuum body as:
\begin{equation}
\begin{split}    
    \int\limits_{V} \delta e_{ij} \sigma_{ij} \dd{V} + \int\limits_{V} \delta \kappa_{i} \mu_{i} \dd{V} =  
    \int\limits_{V} \delta u_i f_i \dd{V} + \\
    \int\limits_{S} \delta u_i t_i \dd{S} + \int\limits_{S} \delta \theta_i m_i \dd{S} - \int\limits_{V} \delta u_i \rho \ddot{u}_i \dd{V}.
\end{split}
    \label{eq:dalembert1}
\end{equation}

In this expression, the internal, external, and inertial virtual work contributions are expressed in terms of the virtual displacement \(\delta u_i\) and virtual rotation \(\delta \theta_i\), along with their associated measures of deformation and stress. These virtual variables are kinematically admissible and vanish on prescribed boundary conditions.

From a mathematical perspective, these virtual displacements align directly with the test functions used in the weak formulation. The variational framework is derived by requiring the governing equations to hold for all \(\delta u_i\) and \(\delta \theta_i\), which are equivalent to test functions in the admissible Sobolev space \cite{kreyszig_introductory_1991}. Thus, the D'Alembert principle not only provides the physical foundation for equilibrium but also serves as the starting point for the functional formulation of the problem.

Rewriting and replacing with the constitutive equations~(\ref{eq:constitutive_aniso}), we get
\begin{equation}
\begin{split}
    \int\limits_{V} \delta e_{ij} C_{ijkl} e_{kl} \dd{V} + \int\limits_{V} \delta \kappa_{i} D_{ij} \kappa_{j} \dd{V} + \int\limits_{V} \delta u_i \rho \ddot{u}_i \dd{V} - \\ \int\limits_{V} \delta u_i f_i \dd{V} - \int\limits_{S} \delta u_i t_i \dd{S} - \int\limits_{S} \delta \theta_i m_i \dd{S} = 0,
\end{split}
    \label{eq:dalembert2}
\end{equation}
where it is required for \(u_i \in H^2(V) \times C^2((0, T])\), due to the definition of \(\kappa_i\) in terms of the displacement, described in section~\ref{subsec:ccst_equations}. To relax this condition, we propose a mixed FEM formulation in terms of the displacement and rotation fields, as presented in previous works \cite{darrall_finite_2014, guarin-zapata_variational_2021}. We assume \(\theta_i\) to be independent of \(u_i\), and enforce their compatibility by adding the following Lagrange multiplier to the equation (\ref{eq:dalembert2})
\begin{equation}
    \int\limits_{V}\lambda_i (\epsilon_{ijk} u_{k,j} - 2\theta_i )\dd{V},
\end{equation}
with the corresponding variation given by
\begin{equation}
    \int\limits_{V}\delta \lambda_i (\epsilon_{ijk} u_{k,j} - 2\theta_i )\dd{V} + \int\limits_{V} (\epsilon_{ijk} \delta u_{k,j} - 2\delta \theta_i ) \lambda_i \dd{V}\, ,
    \label{eq:lagrangevariation}
\end{equation}
which can be added to equation (\ref{eq:dalembert2}) to guarantee field compatibility. This yields
\begin{equation}
\begin{split}
    \int\limits_{V} \delta e_{ij} C_{ijkl} e_{kl} \dd{V} + \int\limits_{V} \delta \kappa_{i} D_{ij} \kappa_{j} \dd{V} + \int\limits_{V} \delta u_i \rho \ddot{u}_i \dd{V} - \\ \int\limits_{V} \delta u_i f_i \dd{V} - \int\limits_{S} \delta u_i t_i \dd{S} - \int\limits_{S} \delta \theta_i m_i \dd{S} + \\ \int\limits_{V}\delta \lambda_i (\epsilon_{ijk} u_{k,j} - 2\theta_i )\dd{V} + \int\limits_{V} (\epsilon_{ijk} \delta u_{k,j} - 2\delta \theta_i ) \lambda_i \dd{V}\ = 0,
\end{split}
    \label{eq:dalembert3}
\end{equation}
which is the weak form we use to derive the FEM method. In this formulation, the Lagrange multiplier equals the skew-symmetric part of the force-stress tensor, and so it is denoted as \(s_i\) from this point onward (\(\lambda_i = s_i\)), for consistency with the literature~\cite{guarin-zapata_variational_2021}. Notice that now it is only required for \(u_i \in H^1(B)\times C^2((0, T]), \theta_i \in H^1(B), s_i \in L^2(B)\), which relaxes the restrictions on these functions compared to previous equations. 

\begin{figure}[t]
\centering
\includegraphics[width= 5 in]{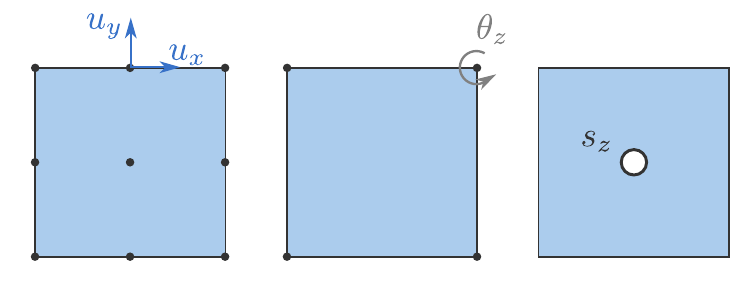}
\caption{Finite element used for the finite element discretization of the
C-CST material model. A second-order Lagrange interpolation is used for
displacements \(\vb{u}\) (left), a first-order Lagrange interpolation for rotations \(\theta_z\) (center), and a piecewise constant is used for the skew-symmetric stresses \(s_z\) (right).}
\label{fig:element}
\end{figure}

To discretize equation~(\ref{eq:dalembert3}), we use second-order Lagrange interpolation functions for the displacement field (\(x, y\) components), first-order interpolation functions for the rotations field, and piecewise constant skew-symmetric stresses. This choice motivated by the need to ensure compatibility between the variables and numerical stability. Specifically, second-order interpolation for displacements captures complex deformations effectively, first-order interpolation for rotations avoids over-constraining the system while maintaining compatibility, and piecewise constant stresses prevent numerical oscillations. Nevertheless, this decision is based on experimental considerations, and a rigorous discussion from a mathematical stability standpoint falls outside the scope of this work. Figure~\ref{fig:element} depicts a typical element for the discretization and the degrees of freedom used.

We follow the combined index notation introduced in \cite{guarin-zapata_variational_2021} to write the discretized equations. Subscripts indicate scalar components of tensors, and capital superscripts indicate interpolation operations. As an example, in the expression
\[u_i = _{u}\hspace{-4pt}N_{i}^{Q}u^{Q},\]
\(i\) indicates the scalar components of the vector \(\vb{u}\), and is added in the shape function to facilitate notation. \(Q\) indicates the nodal point, where \(u^Q\) indicates the corresponding nodal point displacements. The nodal vector implicitly contains both \(x\) and \(y\) Cartesian components, and the interpolation operator is modified accordingly to address this. In \ref{app:interpolation_operators} we present the explicit forms for these operators for further clarification.

The discretized versions of \(u_i, \theta_i\) and \(s_i\) in the described notation are
\begin{equation}
u_i = _{u}\hspace{-4pt}N_{i}^{Q}u^{Q},\quad \theta_i
  = _{\theta}\hspace{-4pt}N_i^{Q}\theta^{Q} ,\quad s_i
  = _{s}\hspace{-4pt}N_{i}^{Q} s^{Q} \, ,
\label{eq:interp_fun_prim}
\end{equation}
and \(e_{ij} , \epsilon_{ijk} u_{i,j}\) and \(\kappa_{i}\) are given by
\begin{equation}
e_{ij} = _{e}\hspace{-4pt}B_{ij}^{Q}u^{Q} ,\quad \epsilon_{ijk} u_{i,j}
  = _{\nabla}\hspace{-4pt}B_{k}^{Q}u^{Q} ,\quad \kappa_{i}
  = _{\kappa}\hspace{-4pt}B_{i}^Q \theta^Q \, ,
\label{eq:interp_fun_second}
\end{equation}

Notice that, according to the finite element presented in Figure~\ref{fig:element}, \(u^Q \in Q_2(\mathcal{Q}), \theta^Q \in Q_1(\mathcal{Q}), s^Q \in Q_0(\mathcal{Q})\), where \(Q_k\) denotes the functional space of quadrilateral elements with \(k\) degree form polynomials. \(\mathcal{Q}\) is a mesh for the domain \(B\) made of quadrilaterals.

Substitution of the above relations in equation~(\ref{eq:dalembert3})
gives the discrete version of the weak form proposed
\begin{equation}
\begin{split}
&\delta u^{Q} \left[\int\limits_{V} (_{e}B_{ij}^{Q}) (C_{ijkl}) (_{e}B_{kl}^{P}) \dd{V}\right] u^{P} + \rho \delta u^{Q} \left[ \int\limits_{V} (_{u}N_{i}^{Q}) (_{u}N_{i}^{P}) \dd{V}\right] \ddot{u}^{P} \\
&- \delta u^{Q}\int\limits_{V}\ _{u}N_{i}^Q f_i \dd{V} - \delta u^{Q}\int\limits_{S}\ _{u}N_{i}^Q t_i \dd{S} 
+ \delta \theta^{Q} \left[\int\limits_{V} (_{\kappa}B_{i}^{Q}) (D_{ij}) (_{\kappa}B_{j}^{P}) \dd{V}\right] \theta^{P}  \\
&-\delta \theta^{Q}\int\limits_{S}\ _{\theta}N_i^Q m_i \dd{S} + \delta s^{Q} \left[\int\limits_{V} (_{s}N_{k}^{Q}) (_{\nabla}B_{k}^{P}) \dd{V}\right]  u^{P} \\
&+ \delta u^{Q} \left[\int\limits_{V} (_{\nabla}B_{k}^{Q}) (_{s}N_{k}^{P})  \dd{V} \right]s^{P} -\delta s^{Q} \left[\int\limits_{V} 2 (_{s}N_{k}^{Q}) (_{\theta}N_{k}^{P}) \dd{V}\right] \theta^{P} \\ &  -\delta \theta^{Q} \left[\int\limits_{V} 2(_{\theta}N_{k}^{Q}) (_{s}N_{k}^{P}) \dd{V}\right] s^{P} = 0 \, .
\end{split}
\label{eq:discrete_PVW}
\end{equation}

From the arbitrariness in the variations \(\delta u^Q\) , \(\delta \theta^Q\)
and \(\delta s^Q\) in equation~\eqref{eq:discrete_PVW} it follows that: 
\begin{equation*}
\begin{split}
&\int\limits_{V} (_{e}B_{ij}^{Q}) (C_{ijkl}) (_{e}B_{kl}^{P}) \dd{V} u^P
  + \rho\int\limits_{V} (_{u}N_{i}^{Q}) (_{u}N_{i}^{P}) \dd{V} \ddot{u}^P \\
 & - \int\limits_{V}\ _{u}N_{i}^Q f_i \dd{V} - \int\limits_{S}\ _{u}N_{i}^Q t_i \dd{S} = 0\,, \\
&\int\limits_{V} (_{\kappa}B_{i}^{Q}) (D_{ij}) (_{\kappa}B_{j}^{P}) \dd{V} \theta^P 
  - \int\limits_{S}\ _{\theta}N_i^Q m_i \dd{S}
  - \int\limits_{V} 2(_{\theta}N_{k}^{Q}) (_{s}N_{k}^{P}) \dd{V} s^P=0\, , \\
&\int\limits_{V} (_{s}N_{k}^{Q}) (_{\nabla}B_{k}^{P}) \dd{V}  u^P
  - \int\limits_{V} 2 (_{s}N_{k}^{Q}) (_{\theta}N_{k}^{P}) \dd{V} \theta^P = 0\, ,
\end{split}
\end{equation*}
which can be written in the standard finite element form for dynamic equilibrium as
\begin{equation}
\begin{bmatrix}
K_{uu}^{QP} &0 &K_{us}^{QP}\\
0 &K_{\theta\theta}^{QP} &-K_{\theta s}^{QP}\\
K_{s u}^{QP} &-K_{s\theta}^{QP} &0
\end{bmatrix}
\begin{Bmatrix}
u^P\\
\theta^P\\
s^P
\end{Bmatrix}
+
\begin{bmatrix}
 M_{uu}^{QP} &0 &0\\
0 &0 &0\\
0 &0 &0
\end{bmatrix}
\begin{Bmatrix}
\ddot{u}^P\\
\theta^P\\
s^P
\end{Bmatrix}
= \begin{Bmatrix}
F_u^Q\\
m_\theta^Q\,\\
0\end{Bmatrix}
\label{eq:mat_fem17}
\end{equation}
where the individual terms are defined in \ref{app:semi-discrete-problem}.

Notice that the \(\ddot{u}\) time derivative is still pending to be discretized up to this point.

\subsection{Time-marching scheme}
\label{subsec:time_marching_scheme}

Equation~(\ref{eq:mat_fem17}) can be rewritten as
\begin{subequations}
    \label{eq:time_system}
    \begin{align}
        [M_{uu}] \{\ddot{u}\} + [K_{uu}]\{u\} + [K_{us}] \{s\} = \{F_u\}\, , \label{eq:time_system1}\\
        [K_{\theta \theta}] \{\theta\} - [K_{\theta s}] \{s\} = \{m_\theta\}\, , \label{eq:time_system2}\\
        [K_{su}] \{u\} - [K_{s \theta}] \{\theta\} = \{0\} , \label{eq:time_system3}
    \end{align}
\end{subequations}

Matrices and vectors in equations~(\ref{eq:time_system}) are obtained after global assembly of the local matrices presented in \ref{app:semi-discrete-problem}. Notice that this system differs significantly from the ones where time-integration methods are typically implemented due to the addition of mixed and Lagrange-multiplier variables in the method used, so an additional procedure is required.

From equation (\ref{eq:time_system2}) one has that \(\{\theta\}\) is
\begin{equation}
    \{\theta\} = [K_{\theta \theta}]^{-1} [\{m_{\theta}\} + [K_{\theta s}] \{s\}]\,
    \label{eq:time_system4}
\end{equation}
and this can be replaced in (\ref{eq:time_system3}) to derive
\begin{equation}
    \{s\} = [C_2]^{-1} [[K_{su}]\{u\} - [C_1]],
    \label{eq:time_system5}
\end{equation}
where $[C_1] = [K_{s \theta}][K_{\theta \theta}]^{-1}\{m_{\theta}\}$ and $[C_2] = [K_{s \theta}][K_{\theta \theta}]^{-1}[K_{\theta s}]$, which are support variables used to keep the expressions readable. Replacing (\ref{eq:time_system5}) in (\ref{eq:time_system1}) we obtain
\begin{equation}
    [M_{uu}]\{\ddot{u}\} + [C_3]\{u\} - [C_4]= 0,
    \label{eq:time_system6}
\end{equation}
where $[C_3] = [K_{uu}] + [K_{us}][C_2]^{-1}[K_{su}]$ and $[C_4] = \{F_u\} + [K_{us}][C_2]^{-1}[C_1]$. This equations is convenient for the implementation of the time-integration scheme. For this purpose, we apply a backward
second-order finite difference approximation for the displacement time derivative
\begin{equation}
\begin{split}
    &\{\ddot{u}\} = \frac{1}{\Delta t^2} \{ \{u^{n+1}\} - 2\{u^n\} + \{ u^{n-1} \} \}, \\
    &\{u\} = \{u^{n+1}\},
\end{split}
    \label{eq:backward_approx}
\end{equation}
where the superindices denote the time step of node values of the displacement field. Replacing the approximation (\ref{eq:backward_approx}) in (\ref{eq:time_system6}) and grouping corresponding terms we obtain
\begin{equation}
\begin{split}
    [ [M_{uu}] + \Delta t^2 [C_3] ] \{u^{n+1}\} = \Delta t^2 [C_4] + [M_{uu}]\{2\{u^n\} - \{u^{n-1}\}\} \,
\end{split}
    \label{eq:time_integration}
\end{equation}
which can be implemented iteratively.

%% file: sections/03_results.tex
\section{Results}
\label{sec:results}

We analyze the simulation results from the computational implementation of the previously described methods. First, we check the stability of the FEM discretization and its implementation through static case experiments. Then, we test the time-marching scheme using dynamic systems. The computational implementations were built on the in-house finite element library SolidsPy \cite{guarin-zapata_nicoguarosolidspy_2021}, and are publicly available at the following GitHub repository: \url{https://github.com/jhortizo/time-domain-ccst}

\subsection{Static verification}

\begin{figure}[t]
\centering
\includegraphics[width=4 in]{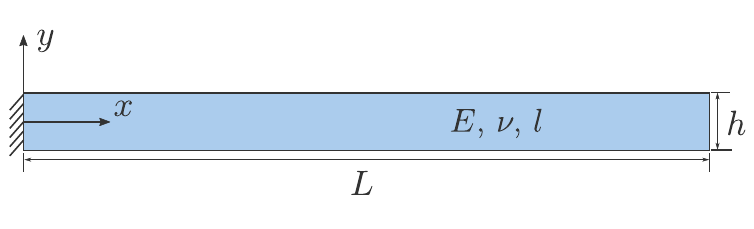}
\caption{Schematic of the cantilever beam. A distributed load is applied on the right side, and the length scale (\(l\)) is varied to verify the behavior of the effective rigidity.
}\label{fig:cantilever-scheme}
\end{figure}

As the first case study, we reconstruct the transverse deformation behavior of a cantilever as the geometric scale approaches the length scale \(l\) (introduced in equation~(\ref{eq:lengthscale})), as presented in \cite{darrall_finite_2014}. Mimicking the reference paper, we set \(E = 2\) and \(\nu = 0\) for a cantilever with height \(h=1\), considered the geometric scale of the system, and lengths \(L = 20h\) and \(L = 40h\) in two evaluated scenarios. Figure~\ref{fig:cantilever-scheme} illustrates the system.

We fixed the \(x, y\) displacements and rotation of the cantilever at the left end. Instead of applying a unit displacement on the cantilever's right side, we imposed a distributed unit load. To achieve different \(h/l\) ratios, we set \(h=1\) and varied the value of \(\eta\) to obtain the required \(l\) values within the desired ratio interval.

In this scenario, the C-CST theory predicts that the effective rigidity of the system \(K\), defined as the ratio between the vertical reaction force and the vertical displacement at the point, increases compared to its classical value of \(\frac{3EI}{L^3}\) as the ratio \(h/l\) increases. This size-dependent behavior arises from the additional effects of couple-stress interaction, which enhances the effective rigidity of the cantilever.

\begin{figure}[t]
\centering
\includegraphics{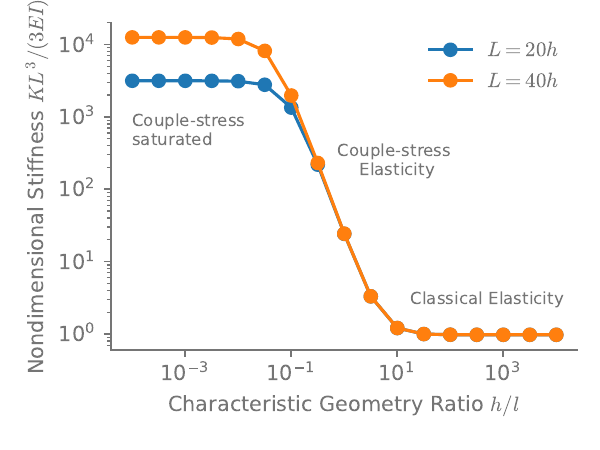}
\caption{Variation of a cantilever effective rigidity for different $h/l$ ratios. This behavior is congruent with previous literature results \cite{darrall_finite_2014}.}\label{fig:cantilever-rigidity}
\end{figure}

Figure~\ref{fig:cantilever-rigidity} presents the obtained results. As previously reported, we can distinguish three domains of interest: a classical domain where the length scale is too small compared to the system's geometric scale, a couple-stress saturated domain where the scenario is reversed, and an intermediate domain where the C-CST model differentially affects the classical case. As expected, \(K = \frac{3EI}{L^3}\) in the classical domain.

The obtained results align with those presented in \cite{darrall_finite_2014}, verifying that the proposed model and its computational implementation agree with previous results.

\subsection{Method of Manufactured Solutions for static case}

We verify the static case using the Method of Manufactured Solutions (MMS)~\cite{aycock_method_2020}. For the MMS, we propose an arbitrary vector function $u$ and apply the differential operator from equation~\ref{eq:ccst-pde} to calculate the body force function. This function is properly discretized and subsequently applied as loads in the implemented FEM method. We then compare the proposed and retrieved solutions, verifying that the error decreases monotonically as the number of mesh elements increases. 

\begin{figure}[t]
\centering
\includegraphics{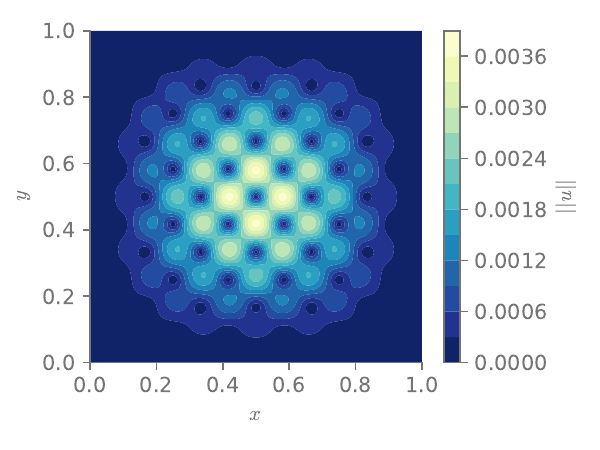}
\caption{Euclidean norm of the manufactured solution proposed.}\label{fig:manufactured-solutions-u}
\end{figure}

We specifically chose to solve for a \(1 \times 1\) rectangle with fixed null displacements and rotations at the boundaries. The proposed solution, which incorporates these boundary conditions, is defined as
\begin{equation}
    \mathbf{u} = \left( x - x^2 \right)^2 \left( y - y^2 \right)^2\begin{Bmatrix}
     \sin(6 \pi x) \cos(6 \pi y) \\
     \cos(6 \pi x) \sin(6 \pi y) \\
    0
    \end{Bmatrix}
    \label{eq:manufactured-solution-u}
\end{equation}
where the terms \(\left( x - x^2 \right)^2\) and \(\left( y - y^2 \right)^2\) ensure that the displacement function cancels at the boundaries, while the sinusoidal functions introduce higher frequency components to assess the convergence of the solver. Figure~\ref{fig:manufactured-solutions-u} presents the Euclidean norm of the displacement field proposed as manufactured solution.

\begin{figure}[t]
\centering
\includegraphics{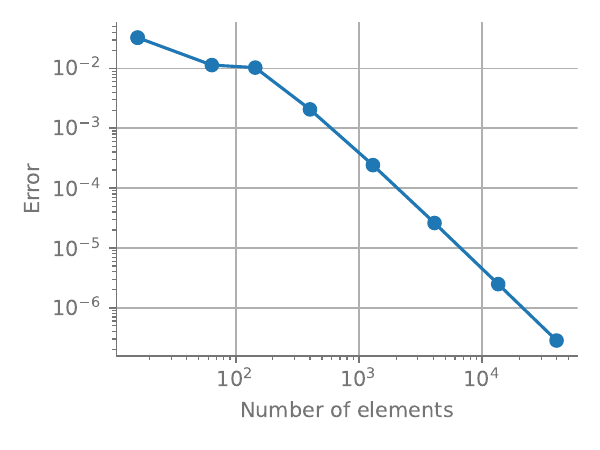}
\caption{\(L_2(V)\) absolute error convergence for a different number of elements, in the Method of Manufactured solutions for the static case. The slope of the linear region is \(-1.88\).}\label{fig:manufactured-solutions}
\end{figure}

The convergence plot for different numbers of elements in the mesh is presented in Figure~\ref{fig:manufactured-solutions}. These numbers of elements correspond to mesh sizes spaced evenly on a log10 base from 0 to $10^{-2}$, approximately 0.56, 0.31, 0.17, 0.1, 0.05, 0.03, 0.017 and 0.01.    

Notice the plot monotonically decreases with a slope of \(-1.88\). This constant monotonicity shows that the method is stable. 

\subsection{Eigenstates' dynamic stability}

As verification for the time-domain solver, we explored the behavior of the eigenstates of the cantilever with support on one side, i.e., fixed \(u_x\), \(u_y\), and \(\theta\) on the left side. We compared several eigenstates with their counterparts from the classical theory, which were also propagated in time using an implicit scheme with second-order backward finite differences for the time derivatives. Corresponding animations for all scenarios are available in the \href{https://figshare.com/articles/media/Supplementary_material_for_i_Implicit_Time-Marching_for_Lagrange_i_i_Multiplier_Formulation_for_Couple_Stress_i_i_Elastodynamics_i_/28281209}{Supplementary Material}.

We simulated a system with height \(h=1\), length \(L=10\), and 48 elements. With \(E=1\), \(\nu=0.29\), \(\rho=1\), and \(\eta=0.1\) (the last considered only in the C-CST case), the simulations were carried out with a time step \(\Delta t = 0.5\) across 1000 iterations. To present dynamic results statically in the paper, we focused on the \(x\), \(y\) components or norm of the displacement at sampling points located along the central axis of the cantilever, equidistant from each other. For reference, the positions of the sample points and the initial state of the system are also shown in the plot.

\begin{figure}[t]
\centering
\includegraphics[width=0.8\linewidth]{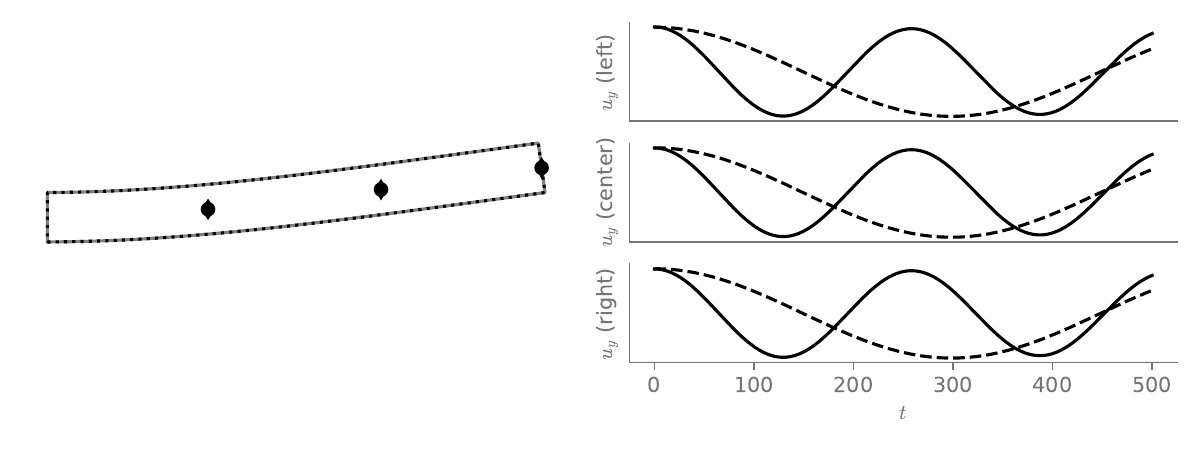}
\caption{Dynamic evolution of cantilever with support on one side, on its first eigenstate. On the left we present the initial state of the simulation, which corresponds to the eigenstate, both for the C-CST (black continuous line) and classical (gray dashed line) models, and with the displacement normalized. On the right, the evolution of $u_y$ over time, without normalization, for both C-CST (continuous lines) and classical (dashed lines) models. The corresponding animation is available in the Supplementary Material.}\label{fig:compare-mode-0}
\end{figure}

We present the first eigenstate dynamic evolution in Figure~\ref{fig:compare-mode-0}. Notice the evolution of the points stays bounded and presents an oscillatory behavior. Even more, the amplitudes and phases of the periodic displacements correspond with the expected behavior. This indicates the method is stable over time and presents simulated results in agreement with theory, which indicates the points, and the whole cantilever, oscillate within the eigenstate. 

Another interesting observation comes from comparing the periodicity of the displacements from C-CST and classical cases. We notice an increase in the natural frequency of the state in the C-CST model, which implies a related increase in the corresponding eigenvalue of the state. From a mechanical point of view, this can be interpreted as a side effect of the additional couple interactions between the material points, just as the increased effective rigidity presented in previous results. The increased rigidity is in agreement with the increase in the natural frequencies of the system, from a mechanical perspective.

\begin{figure}[t]
\centering
\includegraphics[width=0.8\linewidth]{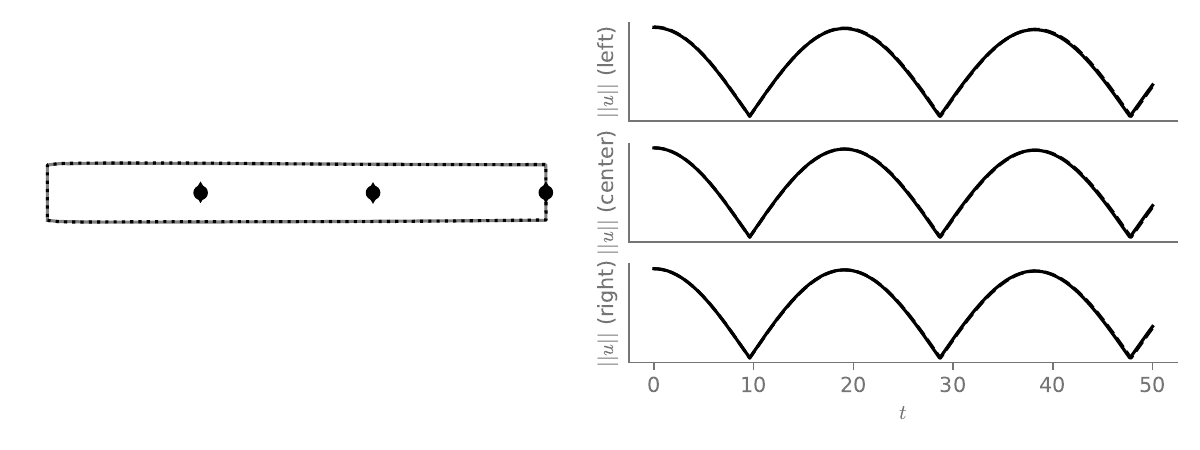}
\caption{Dynamic evolution of cantilever with support on one side, on its third eigenstate (longitudinal oscillation). On the left we present the initial state of the simulation, which corresponds to the eigenstate, both for the C-CST (black continuous line) and classical (gray dashed line) models, and with the displacement normalized. On the right, the evolution of \(\Vert\mathbf{u}\Vert\) over time, without normalization, for both C-CST (continuous lines) and classical (dashed lines) models. The corresponding animation is available in the \href{https://figshare.com/articles/media/Supplementary_material_for_i_Implicit_Time-Marching_for_Lagrange_i_i_Multiplier_Formulation_for_Couple_Stress_i_i_Elastodynamics_i_/28281209}{Supplementary Material}.}\label{fig:compare-mode-2}
\end{figure}

Then, we present the evolution of the third eigenstate in Figure~\ref{fig:compare-mode-2}, using \(\Delta t = 0.05\) for this case and the following. This corresponds to a longitudinal oscillation, hence not presenting a rotational field. For this reason, both the initial state and dynamic evolution for C-CST and classical cases are identical. This is also a symptom for the time integration strategy proposed, which yields the same classical behavior where required. 

\begin{figure}[t]
\centering
\includegraphics[width=0.8\linewidth]{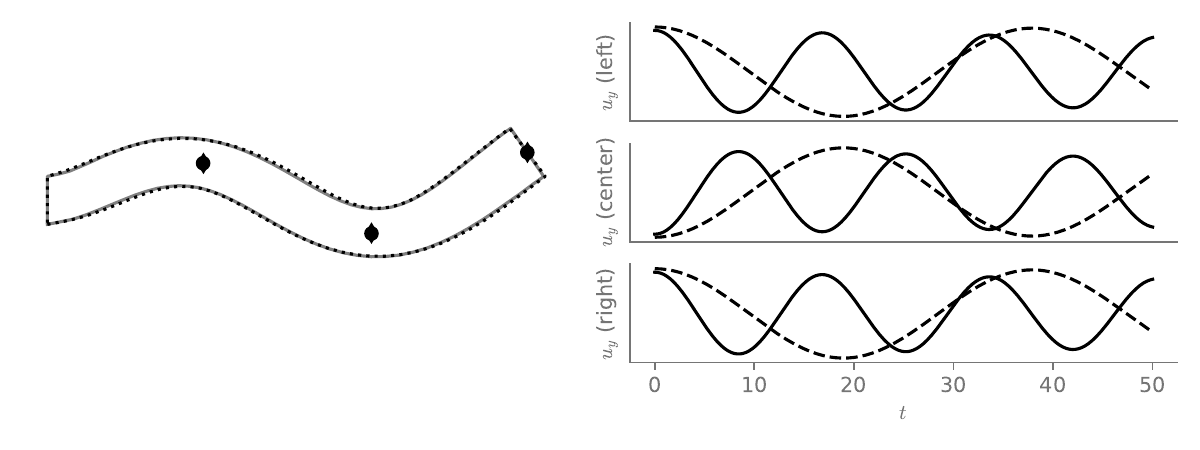}
\caption{Dynamic evolution of cantilever with support on one side, on its fourth eigenstate. On the left we present the initial state of the simulation, which corresponds to the eigenstate, both for the C-CST (black continuous line) and classical (gray dashed line) models, and with the displacement normalized. On the right, the evolution of \(u_y\) over time, without normalization, for both C-CST (continuous lines) and classical (dashed lines) models. The corresponding animation is available in the Supplementary Material.}\label{fig:compare-mode-3}
\end{figure}

Next, in Figure~\ref{fig:compare-mode-3} we present the corresponding results for the fifth eigenstate. As in previous cases, we notice the frequency of the oscillation for the C-CST case increases in contrast to the classical case.

\begin{figure}[t]
\centering
\includegraphics[width=0.8\linewidth]{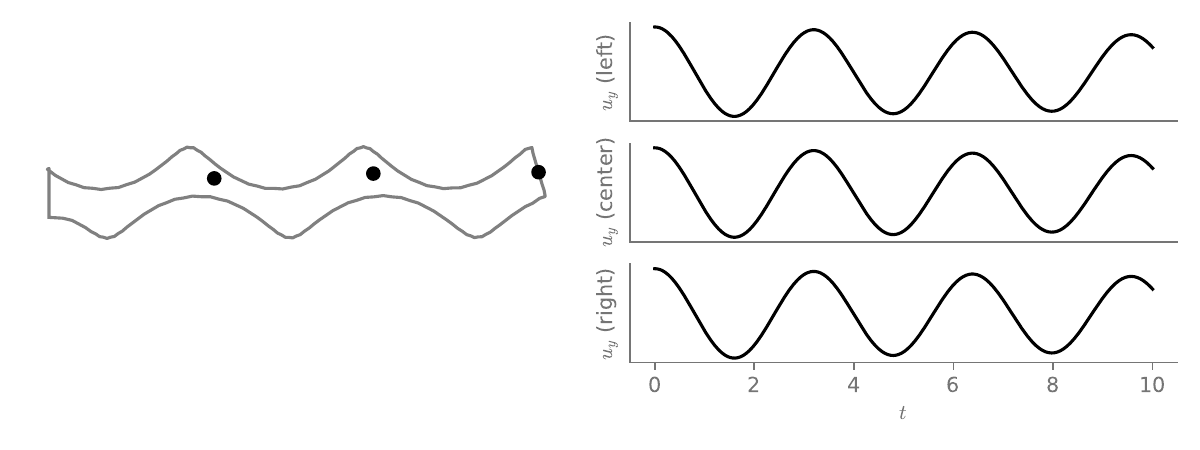}
\caption{Dynamic evolution of cantilever with support on one side, on its twelfth eigenstate. On the left, we present the initial state of the simulation, which corresponds to the eigenstate for the C-CST case with the displacement scaled. On the right, the evolution of $u_y$ over time, without normalization. The corresponding animation is available in the Supplementary Material.}\label{fig:compare-mode-11}
\end{figure}

As a final experiment, we present results for a higher frequency eigenstate, and \(\eta=1, \Delta t = 0.01\), as further proof of the computational stability of the implicit integrator. These results are presented in Figure~\ref{fig:compare-mode-11}.

\subsection{Energy drift analysis}

\begin{figure}[t]
\centering
\includegraphics[width=5 in]{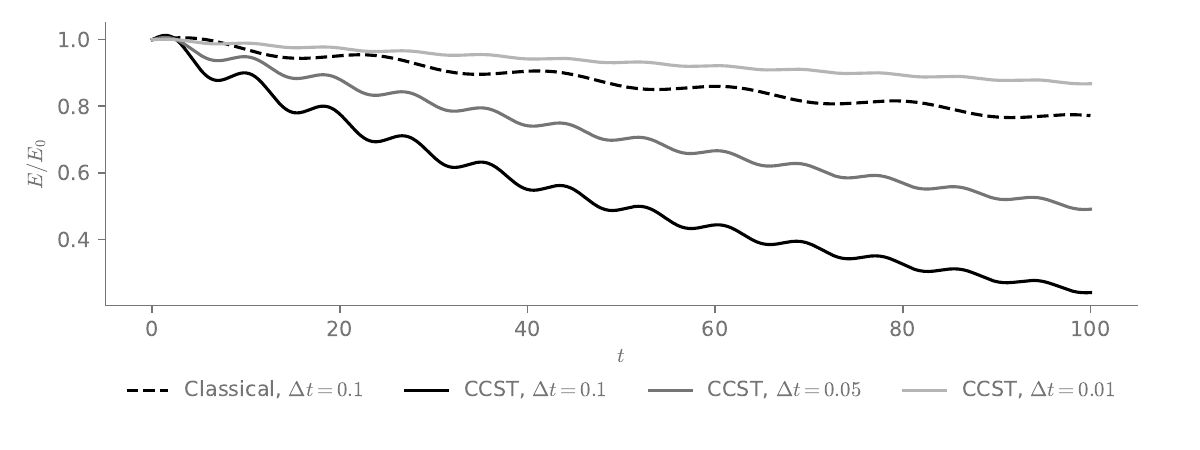}
\caption{Energy evolution for the fifth eigenstate, for both C-CST and Classical models, for different \(\Delta t\). The same case in C-CST loses three times the energy compared to the classical scenario.}\label{fig:energy-comparison}
\end{figure}

Although not evident in the previous figures, the interaction scheme exhibits energy dissipation. In this section, we provide a more detailed analysis of this phenomenon, comparing the dissipation versus \(\Delta t\) and against the same integrator in the classical case. We present the energy evolution for the system discussed in the previous section. This time, we fixed the total simulation time at \(t_f = 100\) and used time step values of \(\Delta t = [0.1, 0.05, 0.01]\). Figure~\ref{fig:energy-comparison} shows the ratio of the energy at each time step relative to the initial state, for the fifth eigenstate presented in Figure~\ref{fig:compare-mode-3}

As shown in Figure~\ref{fig:energy-comparison}, the C-CST integration exhibits more energy drift than its classical counterpart. This is likely due to the additional calculations required for the extra degrees of freedom in the FEM scheme. Consequently, the C-CST case loses, on average, three times more energy per iteration. Notably, even reducing the time step by half does not yield a similar drift, underscoring the increased computational intensity needed to achieve results comparable to the classical case under the current scheme. Only by using one-tenth of \(\Delta t\) could we achieve a smaller drift.

These results emphasize the importance of exploring symplectic integrators for this model, which we leave for future work.

\subsection{Pulse dispersion}

As a final case study, we simulate the propagation of a pulse in C-CST medium consisting of a \(1.5 \times 0.3\) rectangle with 208 elements, to showcase the natural dispersion inherent in the model. Thus, we should see that propagating pulses change shape because different harmonics travel at different phase speeds. Transverse waves propagating in C-CST media exhibit dispersion, with a dispersion relation described by \cite{guarin-zapata_variational_2021}
\begin{equation*}
    \omega = c_2 k \sqrt{1 + k^2 l^2},
    \label{eq:dispersion}
\end{equation*}
where \(\omega\) is the angular frequency, \(k\) is the wavenumber, \(c_2^2 = \mu/\rho\) is the low-frequency phase speed for the transverse wave, and \(l\) is the C-CST length parameter. 

In this case, we propagate a Gaussian pulse that has a spectrum containing components from a range of wavenumbers --- exhibiting the dispersion phenomenon. The domain elements have an average length of 0.05, \(\Delta t = 0.001\), and \(t_f = 3.5\), using material parameters \(E = 1\), \(\nu = 0.29\), \(\eta = 0.001\), and \(\rho = 1\). These material parameters correspond to a length parameter \(l=0.058\), giving a length ratio \(h/l = 0.3/0.058 = 5.1\).

The initial conditions for the displacement and velocity components are:

\begin{equation*}
    \mathbf{u}(x) = \begin{Bmatrix}
     0\\
     e^{-100 \left(x - \frac{L}{2}\right)^2} \\
    \end{Bmatrix},
\end{equation*}

\begin{equation*}
    \mathbf{\dot{u}}(x) = \mathbf{0},
\end{equation*}

while we have the following for the additional degrees of freedom
\begin{align*}
    \theta &= \frac{(100 L - 200 x) e^{- 100 \left(x - \frac{L}{2}\right)^2}}{2}\, ,\\
    s &= 10000 \eta (2x - L ) (50 L^2 - 200 L x + 200 x^2 - 3) e^{- 100 \left(x - \frac{L}{2}\right)^2}\, .
\end{align*}

\begin{figure}[t]
\centering
\includegraphics[width=0.8\linewidth]{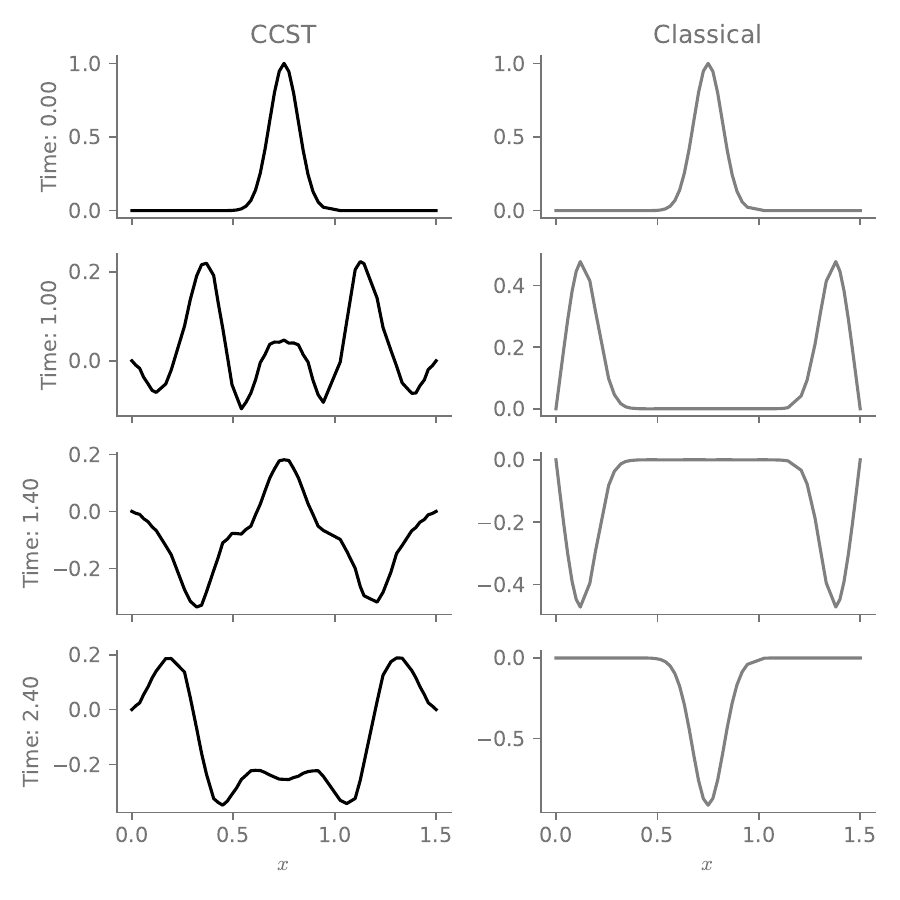}
\caption{\(u_y\) on the cross-section of the rectangle, at \(y=0.15\). Notice that in the C-CST case, the pulse starts losing its form.}\label{fig:pulse-comparison}
\end{figure}

We present four snapshots of the time evolution in Figure~\ref{fig:pulse-comparison}, with the corresponding animation included in the \href{https://figshare.com/articles/media/Supplementary_material_for_i_Implicit_Time-Marching_for_Lagrange_i_i_Multiplier_Formulation_for_Couple_Stress_i_i_Elastodynamics_i_/28281209}{Supplementary Material}.

Figure~\ref{fig:pulse-comparison} illustrates the inherent dispersion of the model, showing how the pulse changes shape as it propagates through the medium for the C-CST case in contrast to the classical one.

%% file: sections/04_conclusions.tex
\section{Conclusions}

The primary contribution of this paper lies in the formulation, implementation, and verification of an implicit time integration scheme using a second-order backward finite difference for obtaining time-domain solutions in C-CST continuum mechanics. This approach eliminates the need for convolutional or frequency-domain calculations. We tested the stability of this scheme through computational simulations of the dynamic evolution of eigenstates in a well-studied physical system, specifically a one-side-supported cantilever. Additionally, physics-oriented simulations confirmed the dispersion of transversal waves during the propagation of a Gaussian pulse in C-CST media. Thus, our proposed method paves the way for further exploration of dynamic effects arising from couple interactions in continuum bodies.

While this work establishes a solid foundation, there remains room for further exploration. A more comprehensive study of the behavior of eigenstates in C-CST systems as the length scale increases is crucial based on our findings. Additionally, implementing more robust time integration schemes, as seen in the current state of the art, could help reduce numerical diffusion and increase the permissible time step size for stable computational simulations. Further verification of the scheme's stability through time-domain MMS is recommended. Finally, additional mechanical case studies, particularly those that were previously restrictive, challenging to simulate, or unfeasible with earlier dynamic methods, should be investigated~\cite{lei_laplace-domain_2023, dargush_convolved_2023, irwin_large_2024, sharma_dynamical_2016, repka_numerical_2019}.

%% file: sections/appendix.tex
\appendix

\section{Interpolation operators used for FEM method}
\label{app:interpolation_operators}
We have the following explicit forms for the interpolation matrices in two
dimensions \cite{bathe_finite_1996}: 
\begin{align*}
_{u}\mathbf{N}^Q &= _{\theta}\mathbf{N}^Q = _{s}\hspace{-2pt}\mathbf{N}^Q =
\begin{bmatrix}
\cdots &N^Q  &  0 &\cdots\\
\cdots & 0 &N^Q &\cdots
\end{bmatrix}\, ,
&_{e}\mathbf{B} = \begin{bmatrix}
\cdots & \frac{\partial N^Q}{\partial x}  & 0 &\cdots\\
\cdots & 0 &\frac{\partial N^Q}{\partial x}  &\cdots\\
\cdots & \frac{\partial N^Q}{\partial y}  & \frac{\partial N^Q}{\partial 
x}&\cdots
\end{bmatrix}\, ,\\
_{\kappa}\mathbf{B} &= \begin{bmatrix}
\cdots & -\frac{\partial N^Q}{\partial y}  &\cdots\\
\cdots & \frac{\partial N^Q}{\partial x}  &\cdots
\end{bmatrix}\, ,
&\mathbf{_{\nabla}B} = \begin{bmatrix}
\cdots & -\frac{\partial N^Q}{\partial y}  &\frac{\partial N^Q}{\partial x} 
&\cdots
\end{bmatrix}\, .
\end{align*}
an the following constitutive tensors in Voigt notation
\begin{align*}
&\mathbf{C} = \frac{E (1 - \nu)}{(1 + \nu) (1 - 2\nu)}\begin{bmatrix}
1 & \frac{\nu}{1 - \nu} &0\\
\frac{\nu}{1 - \nu} & 1 &0\\
0 & 0 &\frac{1 - 2\nu}{2(1 - \nu)}
\end{bmatrix}\, ,\quad 
&\mathbf{D} = 4\eta \begin{bmatrix}
1 & 0\\
0 &1
\end{bmatrix}\, .
\end{align*}

\section{Individual terms in semi-discrete problem}
\label{app:semi-discrete-problem}

The semi-discrete problem can be written as in equation~(\ref{eq:mat_fem17}), where the individual terms of the matrices are given by

\begin{align*}
K_{uu}^{QP} &= \int\limits_{V} (_{e}B_{ij}^{Q}) (C_{ijkl}) (_{e}B_{kl}^{P}) \dd{V}\, ,
&M_{uu}^{QP} &= \rho \int\limits_{V} (_{u}N_{i}^{Q}) (_{u}N_{i}^{p}) \dd{V}\, ,\\
K_{u s}^{QP} &= \int\limits_{V} (_{\nabla}B_{k}^{Q}) (_{s}N_{k}^{P})  \dd{V}\, ,
&F_{u }^{Q} &= \int\limits_{V}\ _{u}N_{i}^Q f_i \dd{V} + \int\limits_{S}\ _{u}N_{i}^Q t_i \dd{S}\, ,\\
K_{\theta \theta}^{QP} &= \int\limits_{V} (_{\kappa}B_{i}^{Q}) (D_{ij}) (_{\kappa}B_{j}^{P}) \dd{V}\, ,
&K_{\theta s}^{QP} &= \int\limits_{V} 2(_{\theta}N_{k}^{Q}) (_{s}N_{k}^{P}) \dd{V}\, ,\\
m_{\theta}^{Q} &= \int\limits_{S}\ _{\theta}N_i^Q m_i \dd{S}\, ,
&K_{s u}^{QP} &= \int\limits_{V} (_{s}N_{k}^{Q}) (_{\nabla}B_{k}^{P}) \dd{V}\, ,\\
K_{s \theta}^{QP} &= \int\limits_{V} 2 (_{s}N_{k}^{Q}) (_{\theta}N_{k}^{P}) \dd{V}\, .
\end{align*}